\documentclass[a4paper]{amsart}
\input{article.sty}

\title{Gradient estimates under integral Ricci bounds}
\date{\today}

\author{Ludovico Marini} 
\address[L. Marini]{Dipartimento di Matematica e Applicazioni,
Università degli Studi di Milano\hyp{}Bicocca, Via R. Cozzi 55, I-20125, Milano}
\email[Corresponding author]{l.marini9@campus.unimib.it}

\author{Stefano Pigola}
\address[S. Pigola]{Dipartimento di Matematica e Applicazioni,
Università degli Studi di Milano\hyp{}Bicocca, Via R. Cozzi 55, I-20125, Milano}
\email{stefano.pigola@unimib.it
}
\author{Giona Veronelli}
\address[G. Veronelli]{Dipartimento di Matematica ed Applicazioni, Università degli Studi di Milano\hyp{}Bicocca, Via R. Cozzi 55, I-20125, Milano}
\email{giona.veronelli@unimib.it}


\begin{document}
\begin{abstract}
In this paper we study $W^{1,p}$ global regularity estimates for solutions of $\Delta u = f$ on Riemannian manifolds.
Under integral (lower) bounds on the Ricci tensor we prove the validity of $L^p$-gradient estimates of the form $|| \nabla u ||_{L^p} \le C (|| u ||_{L^p} + || \Delta u||_{L^p})$. 
We also construct a counterexample which proves that the previously known constant lower bounds on the Ricci curvature are optimal in the pointwise sense. 
The relation between $L^p$-gradient estimates and different notions of Sobolev spaces is also investigated. 
\end{abstract}

\maketitle

\section{Introduction}
\label{sec:introduction}
Let $(M, g)$ be an $n$-dimensional Riemannian manifold. 
Fix $p \in (1, + \infty)$. We say that $(M, g)$ supports an $L^p$\textit{-gradient estimate} if there exists a positive constant $C>0$ such that
\begin{equation}
    \label{eq:Lp gradient}
    \Vert \nabla u \Vert_{L^p(M)} \leq C \left( \Vert u \Vert_{L^p(M)} + \Vert \Delta u \Vert_{L^p(M)}\right), \qquad \forall \varphi \in C^\infty_c(M).
\end{equation}
Here $\Delta$ denotes the (negatively defined) Laplace-Beltrami operator; unless explicitly stated all integrals are computed with respect to the Riemannian volume form $d\mu_g$. 

These first order, global estimates have been extensively studied in recent years. Their motivation lies in the regularity theory of elliptic PDEs, but they have connections with the Riesz transform and the theory of Sobolev spaces on Riemannian manifolds. 
% $L^p$-gradient estimates always hold when the manifold is compact, in the non-compact case, their validity strongly depends on the order of integration $p$ and the geometry.
In the Hilbert case $p = 2$, \eqref{eq:Lp gradient} follows from Young inequality and an integration by parts, and holds on all Riemannian manifolds, \cite[Section 3.1]{P2020}.
Thanks to a celebrated result of Coulhon and Duong, \cite{CD2003}, if $p \in (1, 2]$ boundedness properties of the Riesz transform imply the validity of the $L^p$-gradient estimates on any geodesically complete manifolds. 
An alternative proof of this result using methods of PDEs was also found in \cite[Lemma 1.6]{HMRV2021}.
Conversely, when $p>2$ stronger geometrical assumptions are needed. 
Namely, if the Ricci curvature tensor is bounded from below by some constant then \eqref{eq:Lp gradient} is recovered. 
The result is due to Cheng, Thalmaier and Thompson \cite{CTT2018} and relies on probabilistic arguments although there exists an alternative proof by Pigola (and Meda) which relies on the properties of the Riesz transform, see Theorem 8.2 in \cite{P2020}.

$L^p$-gradient estimates are closely related to a second class of functional inequalities known as $L^p$-Calder\'on-Zygmund inequalities, i.e.,
\begin{equation}
    \label{eq:CZ(p)}
    \|\Hess \varphi \|_{L^p}\le C(\|\varphi\|_{L^p}+\|\Delta\varphi\|_{L^p}), \qquad \forall\,\varphi\in C^\infty_c(M).
\end{equation}
By an integration by parts, the validity of \eqref{eq:CZ(p)} implies the corresponding $L^p$-gradient estimate, see \cite[Corollary 3.11]{GP2015} as well as \cite[Theorem 2]{GPM2019}. 
It should be noted that the converse is generally false: in \cite{MV2021} the first and third author were able to construct an example of a manifold which supports $L^p$-gradient estimates but where the Calder\'on-Zygmund inequalities fail for large $p$.
Conditions which ensure the validity of $L^p$-Calder\'on-Zygmund inequality can be found in \cite{GP2015, IRV2019, IRV2020, P2020, MV2022} or the very recent \cite{BDG2022, CCT2022}

In this paper, we establish $L^p$-gradient estimates under integral Ricci lower bounds, that is, the Ricci curvature is allowed some explosion at $- \infty$ as long as it is controlled in a mean integral sense. 
These bounds appear naturally in some isospectral and geometric variational problems as well as in Ricci and K\"ahler-Ricci flows, \cite{CN2015, S2015, TZ2016, B2017, BZ2017, B2018}.
Under (sufficiently small) integral bounds on the Ricci curvature, several properties of manifolds whose Ricci curvature is uniformly bounded from below are recovered. 
For instance Laplacian and Bishop-Gromov comparisons as well as volume doubling and estimates for the isoperimetric and local Sobolev constants are also proved; see \cite{G1988, Y1988, Y1992, PW1997, PW2001}. 

To be more precise, given $u$ a real measurable function and $B_R(x)$ a geodesic ball of center $x \in M$ and radius $R>0$, we have:
\begin{equation*}
    \Vert u \Vert_{p, B_R(x)} \coloneqq \left( \int_{B_R(x)} |u|^p\right)^{\frac{1}{p}}, \qquad   \Vert u \Vert_{p, B_R(x)}^*\coloneqq \left( \fint_{B_R(x)} |u|^p\right)^{\frac{1}{p}}. 
\end{equation*}
Let $\min \Ric (x)$ be the smallest eigenvalue of $\Ric(x) : T_xM \to T_xM$. For $K \in \R$, we define 
\begin{equation*}
    \rho_K(x) \coloneqq \left( \min \Ric (x) - (n-1)K\right)_-, 
\end{equation*}
where $a_- = (|a| - a)/2$ is the negative part of $a$. 
For $R > 0$ and $1 \leq q< +\infty$, we introduce the following quantities
\begin{equation}
    \label{eq:k(q, R, K)}
    k(x, q, R, K) \coloneqq R^2 \Vert \rho_K \Vert_{q, B_R(x)}^* , \qquad k(q, R, K) \coloneqq \sup_{x \in M} k(x, q, R, K),
\end{equation}
which measure the amount of Ricci curvature lying below $(n-1)K$ in an average $L^q$ sense. 
Our main result is the following: 
\begin{mytheorem}
\label{thm:Lp gradient intro}
Let $(M, g)$ be a complete $n$-dimensional Riemannian manifold and let $n < p < + \infty$. There exist constants $\varepsilon = \varepsilon(p, n, K)> 0$, $C(p, n, K) > 1$ such that if $k(p/2, 1, K) \leq \varepsilon$ for some $K\ge 0$,
then
\begin{equation}
    \label{eq:Lp gradient intro}
     \Vert \nabla u \Vert_{L^p(M)}^p \leq C \left( \Vert u \Vert_{L^p(M)}^p + \Vert \Delta u \Vert_{L^p(M)}^p\right) \quad \forall u \in C^\infty_c(M).
\end{equation}
\end{mytheorem}

Our proof goes as follows.
Thanks to a results of Dai, Wei and Zhang, \cite[Theorem 1.9]{DWZ2018}, integral Ricci bounds at order $p/2$ allow to obtain a local gradient estimate which, properly integrated, produces a local $L^{p}$-gradient estimate. 
Using a covering argument allowed by local volume doubling we then extend the local $L^{p}$-gradient estimate to a global one. 

Since $L^p$-gradient estimates hold on any Riemannian manifold when $p\in(1,2]$, by interpolation the previous Theorem \ref{thm:Lp gradient intro} implies the estimates also for $p<n$.

\begin{mycorollary}
\label{cor:Lp gradient intro interpolation}
Let $(M, g)$ be a complete $n$-dimensional Riemannian manifold and let $n < p_0 < + \infty$. There exists a constant $\varepsilon = \varepsilon(p_0, n, K)> 0$ such that if $k(p_0/2, 1, K) \leq \varepsilon$ for some $K\ge 0$,
then, for every $1 < p \leq p_0$, there exist a constant $C(p, n, K)> 1$ such that 
\begin{equation}
    \label{eq:Lp gradient interpolation}
     \Vert \nabla u \Vert_{L^p(M)}^p \leq C \left( \Vert u \Vert_{L^p(M)}^p + \Vert \Delta u \Vert_{L^p(M)}^p\right) \quad \forall u \in C^\infty_c(M).
\end{equation}
\end{mycorollary}

\begin{remark}
\label{rmk:integral ricci vs ricci bounded from below}
Note that $\rho_K(x) = 0$ if and only if $\Ric(x) \ge K(n-1)g_x$ where the inequality is intended in the sense of quadratic forms. 
In particular if the Ricci curvature is bounded from below $\Ric \ge K(n-1) g$, then $k(q, R, K) = 0$ for all $R> 0$ and $q \in [1, +\infty)$. 
Consequently, \Cref{thm:Lp gradient intro} provides an alternative proof of the result by Cheng, Thalmaier and Thompson, \cite{CTT2018}, using only PDEs methods. 
On the other hand, the integral bounds we assume are in general weaker than the usual pointwise bounds; see Remark \ref{rmk:integral bounds} below.
\end{remark}
\begin{remark}
The proofs of our \Cref{thm:Lp gradient intro}, \Cref{cor:Lp gradient intro interpolation} and of the gradient estimate contained in \Cref{lem:local L^p gradient} rely essentially on the Bochner formula, on the local doubling property, on estimates for the local Sobolev constant and on Moser iteration scheme. All of these ingredients are also available on non-smooth RCD spaces (i.e. infinitesimally Hilbertian metric spaces which satisfy a lower bound on the Ricci curvature in a synthetic sense; see e.g. \cite{A2018} and references therein). It would be interesting to investigate if the same strategy bears results also in this setting. 
\end{remark}

In the second part of this paper, for every $2<p<+\infty$ we construct a complete Riemannian manifold on which the $L^p$-gradient estimate fails. 
While several counterexamples to the validity of the $L^p$-Calder\'on-Zygmund inequalities have been found in recent years, \cite{GP2015, L2020, V2020, MV2021}, in the case of $L^p$-gradient estimates the literature is lacking. 
See Section 9 of \cite{P2020} for an extensive account of the topic.
Using a sequence of conformal deformations on separated balls of the Euclidean plane, we are able to obtain the following:

\begin{mytheorem}
\label{thm:counterex intro}
For any $p>2$ there exists a complete $2$-dimensional Riemannian manifold $(\Sigma, g)$ and a sequence of smooth functions $\lbrace u_k \rbrace \subseteq C^\infty_c(\Sigma)$ such that
\begin{equation*}
     \Vert  u_k \Vert_{L^p(\Sigma)}^p + \Vert \Delta u_k \Vert_{L^p(\Sigma)}^p < 1
\end{equation*}
for all $k \in \N$ but
\begin{equation*}
    \Vert \nabla u_k \Vert_{L^p(\Sigma)}^p \geq k. 
\end{equation*}
In other words, the $L^p$-gradient estimate fails on $\Sigma$. 
\end{mytheorem}

Although the counterexample is constructed in dimension 2, using a trick introduced in \cite{HMRV2021}, it extends to arbitrary dimensions. 
Moreover, given an arbitrary increasing function $\alpha : [0, +\infty) \to [0, +\infty)$ such that $\alpha(t) \to + \infty$ as $t \to +\infty$, it is possible to choose $(\Sigma, g)$ so that 
\begin{equation*}
    \Ric(x) \ge -\alpha(r(x)); 
\end{equation*}
see \Cref{rmk:optimality} below for further details. 
This observation shows that the result of Cheng, Thalmaier and Thompson on $L^p$-gradient estimates under Ricci lower bounds, \cite{CTT2018}, is in fact optimal with respect to pointwise bounds.

As a bypass of our construction we also obtain consequences on the theory of Sobolev spaces.
Unlike the Euclidean setting, on a Riemannian manifold one can define several, generally non-equivalent notions of Sobolev spaces. 
We refer to \cite{V2020, IRV2019, IRV2020, HMRV2021, MV2021, MV2022} for extensive details on the topic. 
Here, we recall the following spaces:
\begin{equation*}
   W^{k, p}(M) = \{ u \in L^p(M) : \nabla^j u \in L^p(M), \quad j = 0, \ldots k\}. 
\end{equation*}
and
\begin{equation*}
H^{2,p}(M)=\left\lbrace u \in L^p:\ \Delta u\in L^p\right\rbrace,
\end{equation*}
endowed with their canonical norms.
The covariant derivatives and the Laplacian are interpreted in the sense of distributions. 
In this language, $L^p$-gradient estimates simply state that the following continuous inclusion $H^{2, p}(M) \subseteq W^{1, p}(M)$ holds true.
While it has been recently observed that the inclusion $W^{2, p}(M) \subsetneq H^{2, p}(M)$ might fail to be an equality, \cite{V2020,MV2021}, the construction of \Cref{sec:counterexample} provides an example where $W^{1,p} \not\subseteq H^{2,p}(M)$.

\begin{acknowledgements}
The authors would like to thank Guofang Wei for the helpful comments and corrections on the first draft of the paper.
\end{acknowledgements}

\section{\texorpdfstring{$L^p$}{L^p}-gradient estimates}
\label{sec:gradient estimates}
This section is devoted to proving \Cref{thm:Lp gradient intro}.
Before we proceed with the actual proof, we point out the following facts, which will be repeatedly used in the following.

\begin{remark}
\label{rmk:controll on k}
As pointed out in Section 2.3 in \cite{PW2001} for the case $K = 0$, smallness of $k(q, R_0, K)$ at a fixed scale $R_0$ implies a control on $k(q, R, K)$ for all scales $R > 0$. 
This is a consequence of a volume comparison result contained in \cite[Lemma 10]{BPS2020}. 
Indeed, if $q > n/2$ there exists $\varepsilon = \varepsilon(n, q, K)>0$ such that if $k(q, R_2, K) < \varepsilon$ then for every $0 < R_1 < R_2$ one has
\begin{equation*}
    k(q, R_1, K) \leq 4 \left(\frac{R_1}{R_2}\right)^2 \left(\frac{v_K(R_2)}{v_K(R_1)}\right)^{\frac{1}{q}} k(q, R_2, K),
\end{equation*}
where $v_K(R)$ is the volume of the geodesic ball of radius $R$ in the $n$-dimensional space form of constant curvature equal to $K$. 
Since $v_K(R_1)\sim R_1^n$ then $k(q, R_1, K) \to 0$ as $R_1 \to 0$, i.e., $k(q, R_1, K)$ can be made arbitrarily small.
See Corollary 13 in \cite{BPS2020}.
%The contrary is also true: there exists $\varepsilon(n, q, K)$ such that if $k(q, R_1, K) < \varepsilon$ then for every $0<R_1<R_2$ we have
%\begin{equation*}
 %   k(q, R_2, K) \leq 4 \left(\frac{R_2}{R_1}\right)^2 \left(\frac{v_K(R_1)}{v_K(R_1/2)}\right)^{\frac{1}{q}} k(q, R_1, K), 
%\end{equation*}
%see Section 2.3 in \cite{PW2001}.

Note also that $k(p, r, K) \leq k(q, r, K)$ whenever $p\leq q$.
\end{remark}

Under the assumption of a small $k(q, 1, K)$, we fist prove a local $L^{2q}$-gradient estimate which is obtained integrating a local gradient estimate proved in \cite{DWZ2018}.

\begin{lemma}
\label{lem:local L^p gradient}
Let $(M, g)$ be a complete $n$-dimensional Riemannian manifold. 
Let $q > n/2$. There exists $\varepsilon =\varepsilon(n, q, K) > 0$, $C(n, q)> 1$ and $0 < R_0 \leq 1$ such that if $k(q, 1, K) \leq \varepsilon$ then 
\begin{equation}
    \label{eq:half harnack}
    \sup_{B_{R/2}(x)} |\nabla u |^2 \leq C R^{-2} \left[(\Vert u \Vert^*_{2, B_R(x)})^2 + (\Vert \Delta u \Vert^*_{2q, B_R(x)})^2 \right],
\end{equation}
for all $0<R\leq R_0$, for all $x \in M$ and for all smooth functions $u$ on $B_1(x)$. 
Moreover, if we let $p = 2q$ there exists a constant $D(n, p) > 0$ such that
\begin{equation}
    \label{eq:local Lp gradient}
    \Vert \nabla u \Vert^p_{p, B_{R/2}(x)} \leq D R^{-p} \left(\Vert u \Vert^p_{p, B_R(x)} + \Vert \Delta u \Vert^p_{p, B_R(x)} \right)
\end{equation}
holds for all $x \in M$, $0< R \leq R_0$ and all smooth functions $u$ on $B_1(x)$. 
\begin{proof}
By Theorem 1.9 in \cite{DWZ2018}, there exists a constant $\varepsilon_0(n, q) > 0$ independent of $R_0$ such that if $k(q, R_0, 0) \leq \varepsilon_0$ then \eqref{eq:half harnack} holds for all $0< R \leq R_0$.
By \Cref{rmk:controll on k} we know that if $k(q, 1, K) \leq \varepsilon$ then $k(q, R, K) \lesssim R^{2-n/q}$ as $R\to 0$ and since $\rho_0(x) \leq \rho_K(x) + (n-1)|K|$ we have 
\begin{equation*}
    k(q, R, 0) \leq k(q, R, K) + (n-1)|K|R^2.
\end{equation*}
Hence, if we take $R_0$ small enough we have $k(q, R_0, 0) \leq \varepsilon_0$ which concludes the first part of the Lemma.
The constant $R_0$ depends on $K, n, \varepsilon$ and $\varepsilon_0$.

From \eqref{eq:half harnack} we have
\begin{equation*}
    \sup_{B_{R/2}(x)} |\nabla u |^{2q} \leq C^q R^{-2q} 2^{q-1}\left[(\Vert u \Vert^*_{2, B_R(x)})^{2q} + (\Vert \Delta u \Vert^*_{2q, B_R(x)})^{2q} \right].
\end{equation*}
By H\"older inequality
\begin{equation*}
    \left(\fint_{B_R(x)} u^2\right)^q \leq \fint_{B_R(x)} u^{2q},
\end{equation*}
hence, 
\begin{equation*}
    \int_{B_{R/2}(x)} |\nabla u|^{2q} \leq C^q R^{-2q} 2^{q-1} \frac{\vol(B_{R/2}(x))}{\vol(B_R(x))}\left(\int_{B_R(x)}  u^{2q} + \int_{B_R(x)}  \Delta u^{2q} \right),
\end{equation*}
which proves \eqref{eq:local Lp gradient}. 
\end{proof}
\end{lemma}

The local $L^p$-gradient estimate \eqref{eq:local Lp gradient} extends to the whole manifold using a uniformly locally finite covering of $M$. 
The existence of such covering is a consequence of local volume doubling which holds under integral Ricci bound; see Lemma 10 and subsequent results in \cite{BPS2020}.
The proof is identical to the case of lower bounded Ricci curvature and is therefore omitted; see Lemma 1.1 of \cite{H1999}. 
\begin{lemma}
\label{lem:covering lemma}
Let $(M, g)$ be a complete $n$-dimensional Riemannian manifold and $q > \frac{n}{2}$. 
There exists $\varepsilon = \varepsilon(n, q, K) > 0$ such that if $k(q, 1, K)\leq \varepsilon$ and $0 < 2R \leq 1$, there exists a sequence $\lbrace x_i\rbrace_{i \in \N} \subset M$ such that
\begin{enumerate}[label=(\roman*)]
    \item $M = \bigcup_{i \in \N} B_{R/2}(x_i)$;
    \item $B_{R/4}(x_i) \cap B_{R/4}(x_j) = \emptyset$ if $i \neq j$;
    \item there exists $N \in \N$ such that every $x \in M$ intersects at most $N$ balls $B_R(x_i)$.
\end{enumerate}
Note that the number $N$ depends on the doubling constant and thus, on $n, q$ and $K$. 

\end{lemma}

We are now ready to prove the global $L^p$-gradient estimate.

\begin{myproof}[of \Cref{thm:Lp gradient intro}]
Let $u \in C^\infty_c(M)$ and $\Omega = \supp(u)$ and let $0 < R \le R_0$ small enough such that $2R \leq 1$.
Here $R_0$ is the radius appearing in \Cref{lem:local L^p gradient}.
By \Cref{lem:covering lemma} and precompactness of $\Omega$, there exist $x_1, \ldots, x_h \in M$ such that 
\begin{enumerate}[label=(\roman*)]
    \item $\Omega \subseteq \bigcup_{i = 1}^h B_{R/2}(x_i)$;
    \item every $x \in \Omega$ intersects at most $N$ balls
    $B_R(x_i)$.
\end{enumerate}
Then
\begin{align*}
    \int_M |\nabla u|^{p} &\leq \sum_{i = 1}^h \int_{B_{R/2}(x_i)} |\nabla u|^{p} \leq D R^{-p}  \sum_{i = 1}^h \left(\int_{B_R(x_i)} |u|^{p} + \int_{B_R(x_i)} |\Delta u|^{p}\right) \\
    &\leq D R^{-p} \int_{M} \sum_{i = 1}^h \mathbbm{1}_{B_R(x_i)} \left(
    |u|^{p} + |\Delta u|^{p}\right) \leq D R^{-p}N \left( \int_{M} |u|^{p} + \int_M |\Delta u|^{p}\right),
\end{align*}
which proves the $L^{p}$ gradient estimate.
% We conclude the proof by an interpolation argument: the validity of the $L^{p_0}$-gradient estimate is equivalent to the boundedness of the gradient as a linear operator $\nabla : H^{2, p_0} \to L^{p_0}$.
% It is known from \cite{CD2003} that $L^p$-gradient estimates hold on any complete Riemannian manifold as long as $1<p\leq 2$, in particular, $\nabla : H^{2, 2}\to L^2$ is a bounded linear operator. 
% It suffices to prove that $L^p$-gradient estimates hold for any $2 < p < p_0$.
% Define $\theta \in (0, 1)$ such that
% \begin{equation*}
%     \frac{1}{p} = \frac{\theta}{p_0} + \frac{1-\theta}{2}.
% \end{equation*}
% Then, since $H^{2, p}$ and $L^p$ are interpolation spaces of exponent $\theta$ between $(L^2, L^{p_0})$ and $(H^{2, 2}, H^{2, p_0})$ (see Theorem 6.4.5 of \cite{BL1976}) then $\nabla: H^{2, p} \to L^p$ is a bounded linear operator, i.e., the $L^p$-gradient estimate holds on $M$
\end{myproof}

Using the fact that $L^p$-gradient estimates always hold on complete Riemannian manifolds if $p \in (1, 2]$, \cite{CD2003}, we fill the gap between 2 and $n$ with an interpolation argument. 
In order to do so, we need the following observations which are interesting in themselves. 

\begin{remark}
\label{rmk:extension to H2p}
If $(M, g)$ is a complete Riemannian manifold supporting an $L^p$-gradient estimate for some $p \in (1, +\infty)$, then \eqref{eq:Lp gradient intro} extends with the same constant to all functions in $H^{2,p}(M)$.
Indeed, if $u \in H^{2,p}(M)$, by a result of Milatovic, \cite[Appendix]{GPM2019}, there exists a sequence $\{u_k\} \subseteq C^\infty_c(M)$ such that $u_k \to u$ with respect to the $H^{2,p}$ norm. 
Applying \eqref{eq:Lp gradient intro} to $u_k$, we deduce that $\nabla u_k$ is Cauchy and thus converges in the space of $L^p$ vector fields. 
Testing $\nabla u_k$ against a smooth and compactly supported vector field and taking the limit shows in fact that $\nabla u_k$ converges in $L^p$ norm to the weak gradient $\nabla u$.
% Denote with $X \in \Gamma_{L^p}(M, TM)$ its limit, take $Y$ a smooth and compactly supported vector field we have that
% \begin{equation*}
%     \int_M \langle \nabla u_k, Y \rangle = - \int_M \Div(Y) u_k.
% \end{equation*}
% Taking the limit on both sides we deduce that $X = \nabla u$. 
\end{remark}

\begin{remark}
\label{rmk:gradient estimates vs operators}
First of all, since the heat semigroup is strongly continuous and a contraction on $L^p(M)$ for all $p \in (1, +\infty)$, \cite[Theorem IV.8]{G2016}, it follows by the Hille-Yosida theorem that 
\begin{equation*}
    || u ||_{L^p} \le || (-\Delta + 1)u ||_{L^p},
\end{equation*}
for all $u \in H^{2,p}(M)$. 
As a consequence, the norms $||u||_{L^p} + ||\Delta u ||_{L^p}$ and $|| (-\Delta + 1)u ||_{L^p}$ are equivalent on $H^{2, p}(M)$. 

Next we define the following operator: if $v \in C^\infty_c(M)$, $v \ge 0$ we know there exists $u \in C^\infty(M) \cap H^{2, p}(M)$ such that
\begin{equation}
    (-\Delta + 1)u = v,
\end{equation}
see for instance \cite[Lemma 6.1]{MV2022}.
We define $Tv = \nabla(-\Delta + 1)^{-1} v = \nabla u$. 
If $v \in C^\infty_c(M)$ takes different signs, then $T$ is extended by linearity.
Moreover, if we have an $L^p$-gradient estimate on $M$, then for all $v \in C^\infty_c(M)$
\begin{equation*}
    ||Tv||_{L^p} = || \nabla u ||_{L^p} \le C || (-\Delta + 1)u ||_{L^p} = || v ||_{L^p}
\end{equation*}
where $u = (-\Delta + 1)^{-1} v \in H^{2, p}(M)$. In other words
\begin{equation*}
    T = \nabla (-\Delta + 1)^{-1} : C^\infty_c(M) \subseteq L^p(M) \to \Gamma_{L^p}(M, TM)
\end{equation*}
is a bounded operator. 
Since $C^\infty_c(M)$ is dense in $L^p(M)$, $T$ extends to a bounded operator defined on $L^p(M)$. In particular, since $C^\infty_c(M)$ is dense in $L^p(M) \cap L^q(M)$ then the extension of $T$ to $L^p(M)$ and of $T$ to $L^q(M)$ coincide on the intersection of the two spaces. 
\end{remark}

\begin{myproof}[of \Cref{cor:Lp gradient intro interpolation}]
We know that the $L^2$-gradient estimates hold on any Riemannian manifolds while the $L^{p_0}$ gradient estimates hold on $M$ as a consequence of \Cref{thm:Lp gradient intro}. 
By \Cref{rmk:gradient estimates vs operators}, $T = \nabla(-\Delta + 1 )^{-1}$ extends to a linear operator 
\begin{equation*}
    T: L^{p_0}(M) \cap L^2(M) \to \Gamma_{L^{p_0}}(M, TM) \cap \Gamma_{L^2}(M, TM) 
\end{equation*}
such that
\begin{equation*}
    || T v ||_{L^{p_0}} \le C_1 || v ||_{L^{p_0}}, \qquad || T v ||_{L^{2}} \le C_2 || v ||_{L^{p2}}
\end{equation*}
for all $v \in L^{p_0}(M) \cap L^2(M)$. 
For $\theta \in (0, 1)$ such that
\begin{equation*}
    \frac{1}{p} = \frac{\theta}{p_0} + \frac{1-\theta}{2}, 
\end{equation*}
by the Riesz-Thorin theorem $T$ extends to a bounded linear operator on $L^p(M)$. Hence, the $L^p$-gradient estimates hold on $M$. 
\end{myproof}

\begin{remark}\label{rmk:integral bounds}
As alluded to in the introduction, the integral curvature bounds assumed here are weaker than the classical pointwise bounds. 
An easy example of a Riemannian manifold $(M,g)$ satisfying $\inf_M \min\Ric = -\infty$ but with $k(p,1,0)$ arbitrarily small, can be constructed as follows. 
We let $M=\R^2$ endowed with the conformally flat metric $g=e^{2\varphi}dx^2$, where $\varphi$ is a smooth nonpositive function. 
In the following the sub/superscript $e$ denotes the objects taken with respect to the Euclidean metric. 
In particular $\vol_g(K)\le \vol_e(K)$ for any measurable set $K\subset \R^2$, and $B^g_R(w)\supseteq B^e_R(w)$ for any $R>0$ and $w\in \R^2$. 
Suppose now that $\supp \varphi \in \cup_{n\in \mathbb N}B^e_{1/2}((4n,0))$. 
This guarantees that $B_1^g(w)\subseteq B_2^e(w)$ for any $w\in\R^2$. 
Moreover, given $w\in \R^2$, let $n_w$ be the unique integer (if any) such that $B_{1/2}^e((4n_w,0))$ intersects $B_1^e(w)$. 
Then 
\begin{equation}
\label{eq:vol lower}
\vol_g B_1^g(w)\ge \vol_g B_1^e(w) \ge \vol_g (B_1^e(w)\setminus B_{1/2}^e(4n_w,0))=\frac 34 \pi.
\end{equation}
Fix $a\in(2-\frac 2p,2)$, and $\phi_0\in C^\infty_c(B^e_{1/2}(0,0))$, we define $\varphi(x,y) = \sum_{n\in \mathbb N} \phi_n(x,y)$, where $\phi_n(x,y)=n^{-a}\phi_0(n(x-4n,y))$ if $n \ge 1$. 
On the one hand, since $\Delta_e\phi_0$ attains positive values and since $\Delta_e\phi_n(x,y)=n^{2-a}\Delta_e\phi_0(n(x-4n,y))$, we have that $\Ric_g=-2\Delta_e\varphi$ is lower unbounded. 
On the other hand, 
\begin{align*}
\int_{B_1^g(w)}((\min\Ric)_-)^p d\mu_g = 2^p\int_{B_1^g(w)}((\Delta_e\varphi)_+)^p d\mu_g &\le 2^p \int_{B_2^e(w)}((\Delta_e\phi_{n_w})_+)^p dx^2\\ &= 2^p n_w^{2p-pa-2} \int_{B_{1}^e(0,0)}((\Delta_e\phi_{0})_+)^p dx^2 \\
&\le 2^p \int_{B_{1}^e(0,0)}((\Delta_e\phi_{0})_+)^p dx^2,
\end{align*}
which is uniformly bounded independently from $w$. 
Moreover, choosing an appropriate $\phi_0$, we can assume that the right hand side of the estimate above is arbitrarily small. 
Together with the uniform volume lower bound \eqref{eq:vol lower}, this proves that $k(p,1,0)<+\infty$ and can be made arbitrarily small. 
\end{remark}

\section{A counterexample}
\label{sec:counterexample}
We provide, for any $n\ge 2$ and any $p>2$, an example of a complete Riemannian manifold which does not support any $L^p$-gradient estimate. 

Take $(\Sigma, g ) = (\R^2, \lambda^2 dx^2)$ where $dx^2$ is the usual Euclidean metric on $\R^2$ and $\lambda \in C^\infty(\Sigma)$ such that $0 < \lambda \leq 1$. 
As above, In the following we denote with $\Delta$ and $\nabla$ the Laplace-Beltrami operator and gradient with respect to the metric $g$ while we use $\Delta_e$ and $\nabla_e$ to denote the Euclidean differential operators. 
The spaces $L^p(\Sigma)$ are defined in terms of the Riemannian volume form $d\mu_g$ while $L^p(\R^2)$ are the spaces with respect to the Lebesgue measure $dx^2$.  
Moreover, take $\lambda(x) = 1$ for all $x \in \Sigma \setminus \bigcup_{m \in \N} B_{1/8}(x_m)$ where $\{x_m\} \subseteq \Sigma$ is a sequence of points going at $\infty$ such that their pairwise distance is lower bounded by $1$. 
Since $(\Sigma, g)$ is isometric to $(\R^2, dx^2)$ outside of a countable union of bounded sets whose pairwise distance is uniformly lower bounded, it is a complete Riemannian manifold.
Next, take $\varphi_0 \in C^\infty_c(\Sigma)$ such that
\begin{equation*}
    \begin{cases}
    \varphi_0(x,y) = x+1 \text{ on } B_{1/4}(x_0) \\
    \supp (\varphi_0) \Subset B_{1/2}(x_0)
    \end{cases}
\end{equation*}
and let $\varphi_m(x, y) = \varphi_0(x-m, y)$, for all $m \in \N$.
Then, we take $u_k \in C^\infty_c(\Sigma)$ defined as 
\begin{equation*}
    u_k = \sum_{m = 0}^k \frac{\varphi_m}{2^m}
\end{equation*} 
and note that, since $\lambda \leq 1$, then 
\begin{equation*}
    \Vert u_k \Vert_{L^p(\Sigma)}^p = \int_\Sigma  \sum_{m = 0}^k \frac{|\varphi_m|^p}{2^{mp}} \lambda^2 dx \leq  \sum_{m = 0}^k \frac{\Vert \varphi_m \Vert_{L^p(\R^2)}^p}{2^{pm}} \leq \Vert \varphi_0 \Vert_{L^p(\R^2)}^p \sum_{m = 0}^{+\infty} 2^{-mp} < + \infty,
\end{equation*}
where the estimate is independent of $k$.
Similarly, we have
\begin{align*}
     \Vert \Delta u_k \Vert_{L^p(\Sigma)}^p &= \int_\Sigma  \sum_{m = 0}^k \frac{|\Delta \varphi_m|^p}{2^{mp}} \lambda^2 dx \\
     &=\sum_{m = 0}^k 2^{-mp}\left(\int_{B_{\frac{1}{4}}(x_m)} |\Delta \varphi_m|^p \lambda^2 dx +\int_{\Sigma \setminus B_{\frac{1}{4}}(x_m)}  |\Delta \varphi_m|^p \lambda^2 dx\right). 
\end{align*}
Since $\Delta \varphi_m = \lambda^{-2} \Delta_e \varphi_m$ then
\begin{align*}
    \Vert \Delta u_k \Vert_{L^p(\Sigma)}^p &= \sum_{m = 0}^k 2^{-mp}\left(\int_{B_{\frac{1}{4}}(x_0)} |\Delta_e \varphi_0|^p \lambda^{2-2p} dx +\int_{\Sigma \setminus B_{\frac{1}{4}}(x_0)} |\Delta_e \varphi_0|^p  dx\right) \\
    &\leq \Vert\Delta_e\varphi_0 \Vert_{L^p(\R^2)}^p\sum_{m = 0}^{+\infty} 2^{-mp} < + \infty.
\end{align*}
Once again, the estimate is independent of $k$. 
Now, let $\lambda_\infty(x) = |x|^{2\beta}$ in $B_\delta(x_0)$ for some $0 \leq \delta \ll 1/8$ and note that
\begin{equation*} \int_{B_\delta(x_0)} |\nabla_e \varphi_0|_e^p \lambda_\infty^{2-p} dx = 2\pi \int_0^\delta r^{1-2\beta(p-2)} dr = +\infty,
\end{equation*}
as long as $\beta(p-2) > 1$. 
Here $|x| = r$ denotes the Euclidean distance from the origin. 
Then for any $m \in \N$ we can find some $\varepsilon_m > 0$, $\varepsilon_m \to 0$ as $m \to +\infty$ such that
% letting $\lambda_m(x) = |x|^{2\beta} +  \varepsilon_m $ for all $x \in B_\delta(x_0)$ we have
\begin{equation*}
 \int_{B_\delta(x_0)} |\nabla_e \varphi_0|^p (|x|^2 + \varepsilon_m)^{(2-p)\beta} dx \geq  2^{mp}.
 %    \int_{B_\delta(x_0)} |\nabla \varphi_0|^p (|x|^2 + \varepsilon_m)^{2\beta} dx \geq  2^{mp}.
\end{equation*}
For $x \in B_{1/8}(x_0)$ and $\varepsilon \in [0, 1]$ we define the family of functions $\lambda_\varepsilon \in C^\infty(B_{1/8}(x_0)) $ such that 
\begin{equation*}
    \begin{cases}
    0 < \lambda_\varepsilon \le 1 \\
    \lambda_\varepsilon(x) = (|x|^{2} + \varepsilon)^{\beta} \text{ if } x \in B_{\delta}(x_0) \\
    \supp(1 - \lambda_\varepsilon) \subseteq B_{1/8}(x_0).
    %\\
    %\min_{x \in B_{1/8}(x_0)} \lambda_\varepsilon (x) = \varepsilon^\beta
    \end{cases}
\end{equation*}
%Moreover, we also require that $\lambda_\varepsilon(x)$ is continuous with respect to $\varepsilon$ for all $\varepsilon \in [0, 1]$. 
%Note that the request on the minimum and the continuity with respect to $\varepsilon$ are not strictly necessary in the counterexample but are used to give a pointwise estimate to the Ricci curvature. 
If we let $\lambda \in C^\infty(\Sigma)$ such that
\begin{equation*}
    \begin{cases}
    0 < \lambda \leq 1 \\
    \lambda(x) = 1 \text{ if } x \in \Sigma \setminus \bigcup_{m \in \N} B_{1/8}(x_m) \\
    \lambda(x) = \lambda_{\varepsilon_m}(x-x_m) \text{ if } x \in B_{\delta}(x_m)
    \end{cases}
\end{equation*}
then
\begin{align*}
    \Vert \nabla u_k \Vert_{L^p(\Sigma)}^p &= \int_\Sigma  \sum_{m = 0}^k \frac{|\nabla \varphi_m|^p}{2^{mp}} \lambda^2 dx \geq \sum_{m = 0}^k 2^{-mp} \int_{B_\delta(x_0)} |\nabla \varphi_o|^p \lambda_m^2 dx \\
    &= \sum_{m = 0}^k 2^{-mp} \int_{B_\delta(x_0)} |\nabla_e \varphi_0|^p (|x|^2 + \varepsilon_m)^{(2-p)\beta} dx \geq k. 
\end{align*}
Finally, note that we can suppose without loss of generality that 
\begin{equation*}
     \Vert  u_k \Vert_{L^p(\Sigma)}^p + \Vert \Delta u_k \Vert_{L^p(\Sigma)}^p < 1
\end{equation*}
(if not, rescale $\varphi_0$ and replace the sequence $\{u_k\}$ with a suitable subsequence). 
This concludes the proof of \Cref{thm:counterex intro}.

Observe that if we take $$u_\infty = \sum_{m = 0}^{+ \infty} \frac{\varphi_m}{2^m},$$ then $u_\infty, \Delta u_\infty \in L^p(\Sigma)$ while $|\nabla u_\infty| \not\in L^p(\Sigma)$. In particular $u_\infty \in H^{2, p}(\Sigma)$ while $u_\infty \not \in W^{1, p}(\Sigma)$. 

\begin{remark}\label{rmk:optimality}
%Note that the Ricci curvature of $(\Sigma, g)$ at $x \in\Sigma$ is given by
%\begin{equation*}
%    \Ric(x) = - \Delta_e \log{\lambda} \ dx^2 = -\frac{ \Delta_e \log{\lambda}}{\lambda^2} g = \left(-\frac{ \Delta_e \lambda }{\lambda^3} + \frac{|\nabla_e \lambda|^2}{\lambda^4}\right) g. 
%\end{equation*}
%In particular, if $x \in \Sigma \setminus \bigcup_{m = 0}^{+ \infty} B_{1/8}(x_m)$ then $\Ric(x) = 0$. 
%On the other hand if $x \in B_{1/8}(x_m)$ we have
%\begin{equation*}
 %   \Ric(x) \ge - \frac{\Delta_e \lambda}{\lambda^2}g.
%\end{equation*}
%Since $\lambda_\varepsilon(x)$ is continuous on $B_{1/8}(x_0) \times [0, 1]$, we have $\Delta_e \lambda_{\varepsilon_m} \le C$ where the estimate is uniform in $m$. 
%Moreover, $\lambda_{\varepsilon_m}(x) \ge \varepsilon_m^\beta$ hence
%\begin{equation*}
 %   \Ric(x) \ge -\frac{C}{\varepsilon_m^{3\beta}} g. 
%\end{equation*}
%Note that as $x \to + \infty$, the lower bound diverges at $- \infty$. 
We observe that the choice of the sequence$\{x_m\}$ is quite arbitrary.
In particular, let $\alpha: [0, + \infty) \to [0, + \infty)$ be an arbitrary increasing function such that $\alpha(t) \to \infty$ as $t \to +\infty$. 
If we choose $x_m$ which diverges quick enough to infinity
%such that $C \varepsilon_m^{-3\beta} \leq \alpha(r(x))$ for all $x \in B_{1/8}(x_m)$, 
we can make the lower bound on Ricci arbitrarily small so that 
\begin{equation*}
    \Ric(x) \geq - \alpha(r(x)). 
\end{equation*}
This shows that the result by Cheng, Thalmaier and Thompson, \cite{CTT2018}, is in fact optimal with respect to pointwise lower bounds. 
\end{remark}
Using a simple trick introduced in \cite{HMRV2021}, this counterexample in dimension $2$ can be used to construct counterexamples to the $L^p$-gradient estimate in arbitrary dimensions $n \geq 2$. 
Indeed, take $(\Sigma, g)$ the above Riemannian manifold and $(N, h)$ any $n-2$ dimensional closed Riemannian manifold. 
Consider $ M = \Sigma \times N$ the product manifold and let $\lbrace v_k \rbrace \subseteq C^\infty_c(M)$ be the sequence defined as follows: for $x \in \Sigma$ and $y \in N$, $v_k(x, y) = u_k(x)$ then we have the following corollary. 

\begin{corollary}
\label{cor:counterexample to L^p grad}
For any $n \geq 2$ and $p>2$, there exists a complete $n$-dimensional Riemannian manifold $M$ such that the $L^p$-gradient estimate fails. More precisely, if $(\Sigma, g)$ is the $2$-dimensional manifold constructed above and $(N, h)$ is an $n-2$-dimensional closed Riemannian manifold one can take $M = \Sigma \times N$. 
\end{corollary}

Similarly, if we take $v_\infty(x, y) = u_\infty(x)$ we have the following result. 

\begin{corollary}
\label{cor: consequence on sobolev spaces}
For any $n \ge 2$ and $p > 2$, there exists a Riemannian manifold $M$ and a function $v_\infty \in H^{2, p}(M)$ such that $v_\infty \not\in W^{1, p}(M)$. 
\end{corollary}

%BIBLIOGRAFIA
\bibliography{biblio}
\bibliographystyle{acm}

\end{document}